\documentclass[12pt]{amsart}
\usepackage{amssymb, enumerate}

\theoremstyle{plain}
\newtheorem{thm}{Theorem}[section]
\newtheorem{lemma}[thm]{Lemma}

\theoremstyle{definition}

\theoremstyle{remark}
\newtheorem{remark}[thm]{Remark}

\newenvironment{enumeratea}
{\begin{enumerate}[\upshape (a)]}
{\end{enumerate}}

\newcommand{\G}{{\mathbf G}}

\newcommand{\cI}{{\mathcal I}}
\newcommand{\eqto}{\stackrel{\lower1.5pt\hbox{$\scriptstyle\sim\,$}}\to}

\DeclareMathOperator{\chara}{char}
\DeclareMathOperator{\spec}{Spec}
\DeclareMathOperator{\pic}{Pic}

\begin{document}

\title[Coverings of Deligne--Mumford stacks]{
   On coverings of Deligne--Mumford stacks and
   surjectivity of the Brauer map
}
\author[A. Kresch]{Andrew Kresch}
\address{
   Department of Mathematics,
   University of Pennsylvania,
   Philadelphia, PA 19104
}
\email{kresch@math.upenn.edu}
\author[A. Vistoli]{Angelo Vistoli}
\address{
   Dipartimento di Matematica,
   Universit\`a di Bologna,
   40126 Bologna, Italy
}
\email{vistoli@dm.unibo.it}
\date{January 22, 2003}
\thanks{Research of the first author supported in part by the Wolfgang
Paul program of the Humboldt Foundation; research of the second author
supported in part by the University of Bologna, funds for selected
research topics.}

\maketitle

\section{Introduction} \label{s.intro}

This paper is concerned with finite covers of Deligne--Mumford stacks
by schemes, in connection with the theory of Brauer group. The reader
is referred to \cite{EHKV} for basic references on
algebraic stacks and Brauer groups.
We are primarily concerned with Deligne--Mumford stacks; every
Deligne--Mumford stack, separated and of finite type over a field,
possesses a coarse moduli algebraic space (\cite{KM}).

Mumford, in \cite{Mu}, constructed intersection theory on the
Deligne--Mumford stack $\overline{\mathcal{M}}_g$ of stable curves of
genus~$g$. One tool that he used was the existence of a finite flat
morphism $Z \to \overline{\mathcal{M}}_g$, where $Z$ is a Cohen--Macaulay
scheme of abelian level structures; his construction would have been made
much simpler by the existence of such a morphism with $Z$ smooth.

Afterwards Looijenga, in \cite{Lo}, showed the existence of a finite
morphism $Z \to \overline{\mathcal{M}}_g$, where $Z$ is a smooth projective
scheme, over the field of complex numbers. A more algebraic construction,
working over more general bases, was given by Pikaart and De Jong in
\cite{PDJ}, and another in \cite{ACV}.

In this paper we show that in fact the existence of such a $Z$ is a very
general phenomenon. We prove that a separated Deligne--Mumford stack with
quasi-projective moduli space which is a quotient stack always possesses a
finite flat cover by a scheme (Theorem~\ref{thm.fifl}); if the stack is
smooth, then the scheme cover can also be chosen to be smooth. The
hypothesis of being a global quotient stack with quasi-projective moduli
space is often verified in practice (for example, for
$\overline{\mathcal{M}}_g$, or for stacks of stable maps into a projective
variety).

We also apply this result to the question of the surjectivity of the map
from the usual Brauer group to the cohomological Brauer group of a scheme.
We recall briefly some essential facts from \cite{Gr}. The (Azumaya)
Brauer group of a scheme $X$ is the group of classes of sheaves Azumaya
algebras on $X$, and this group maps, via the Brauer map, into the cohomological
Brauer group, that is, the torsion subgroup of the \'etale cohomology
group $\mathrm{H}^2(X,\mathbb{G}_m)$. An algebraic stack is called a
quotient stack if it is isomorphic to the stack quotient of an algebraic
space by a linear algebraic group scheme.

The connection between stacks and Brauer groups is one of the central
themes of \cite{EHKV}: the Brauer map is surjective (hence an isomorphism)
if and only if some associated algebraic stacks are all quotient stacks.
In op.\ cit.\ it was shown that knowing the Brauer map to be an
isomorphism implies the quotient stack property for tame Deligne--Mumford
stacks that are gerbes, that is, have nonvarying stabilizer group, of
order prime to the characteristic of the ground field.
Our Theorem~\ref{thm.fifl} is used to show (Theorem~\ref{thm.tfae}) that
the surjectivity of the Brauer map implies the quotient stack property for
smooth Deligne--Mumford stacks of finite type over a field, with some
hypothesis on the characteristic. A final, unconditional, result
(Theorem~\ref{thm.affinemod}) guarantees the existence of Zariski
coverings by quotient stacks for quite general classes of smooth
Deligne--Mumford stacks.

\section{Results} \label{s.results}

\begin{thm} \label{thm.fifl}
Let $X$ be a Deligne--Mumford stack, separated and of finite type over a
field $k$, with quasi-projective coarse moduli space. Assume $X$ is a
quotient stack. Then there exists a quasi-projective scheme $Z$ and a
finite flat local complete intersection morphism $Z\to X$, such that the
singular locus of $Z$ is the inverse image of the singular locus of $X$.
\end{thm}

Recall that a Deligne--Mumford stack $X$ over a field $k$ is called
\emph{tame} if the order of the stabilizer group at any geometric point of
$X$ is relatively prime to the characteristic of $k$ (or if $\chara k=0$).
We say that a Deligne--Mumford stack is \emph{generically tame} if it has a
tame dense open substack.

\begin{thm} \label{thm.tfae}
Given a field $k$ and a positive integer $n$, the following two conditions
are equivalent.

\begin{enumeratea}
\item Every smooth separated generically tame Deligne--Mumford stack over
$k$ of dimension $n$ with quasi-projective moduli space is a quotient
stack.

\item The Azumaya Brauer group of any smooth quasi-projective scheme
over $k$ of dimension $n$ coincides with the cohomological Brauer group.
\end{enumeratea}
\end{thm}

\begin{remark}
In any discussion of quotient stacks, it is worth drawing a comparison with
a related condition for algebraic stacks, known as the resolution
property, which asserts that every coherent sheaf admits a surjection from
a locally free coherent sheaf. The resolution property is discussed in
some detail in a recent paper by Totaro \cite{To}. For the stacks in
Theorem~\ref{thm.tfae}, it is known (cf.~\cite{EHKV}, Remarks 2.15 and
4.3) that the condition to be a quotient stack is equivalent to the
resolution property.
\end{remark}

For $n\le 2$ the assertions of Theorem~\ref{thm.tfae} have long been
known to be true.
For general $n$,
statement (b) is given as an open question in \cite{Mi}, and it is widely
conjectured that the statement is true. Techniques developed by Gabber
may soon settle this question; already an announcement of proof by Gabber
has been made (see, e.g., \cite[p.~19]{Co}).
Gabber has settled the case of affine schemes and separated unions of two
affines \cite{Ga}.
So, we have:

\begin{thm} \label{thm.affinemod}
Every smooth separated generically tame Deligne--Mumford stack over a
field, whose moduli space is either affine or the union of
two affine schemes, is a quotient stack.
\end{thm}

\section{Proofs of results}

The proof of Theorem~\ref{thm.fifl} relies upon the following lemma.

\begin{lemma} \label{lem.bertini}
Let $f\colon U\to V$ be a proper morphism of quasi-projective schemes over
an infinite base field $k$, with constant fiber dimension $r>0$. Choose a
projective embedding $U\to \mathbb{P}^N$ for some $N$. Then for
sufficiently large $d$ the intersection of $U$ with a generic hypersurface
of degree $d$ in $\mathbb{P}^N$ is a Cartier divisor in $U$, of
constant fiber dimension $(r-1)$ over $V$ and whose singular locus
is the intersection of the hypersurface with the singular locus of $U$.
\end{lemma}

\begin{proof}
The irreducible components $E\subset f^{-1}(v)$
of fibers of $f$ over geometric points $v\colon \spec \overline{k} \to V$
vary over a bounded family, hence only finite many Hilbert polynomials
occur.
Choosing $d$ sufficiently large, we may assume the sheaf of ideals
$\cI_E$ of any component $E$ of a geometric fiber of $f$
satisfies $\mathrm{H}^i\bigl(\mathbb{P}^N,\cI_E(d)\bigr)=0$ for $i>0$.
Then the codimension of the space of polynomials of degree $d$ vanishing
along some component of a fiber is bounded below by a polynomial
of degree $r$ in $d$.
Since $r>0$, this codimension is positive for sufficiently large $d$.
Hence a generic hypersurface of degree $d$ will not contain any component
of any geometric fiber of $f$. This establishes the assertion concerning
fiber dimension of the intersection.
The generic hypersurface avoids all the associated points of $U$,
meaning the intersection is a Cartier divisor in $U$.
A standard Bertini-type argument establishes the assertion concerning
the singular locus.
\end{proof}

\begin{proof}[Proof of Theorem~\ref{thm.fifl}]
We are easily reduced to the case of an infinite base field $k$.
Since $X$ is a quotient stack, there exists a projectivized
vector bundle $P_0\to X$ and a
representable dense open substack $Q_0\subset P_0$
with $Q_0\to X$ surjective
(see \cite{EG}). Call $d$ the fiber dimension of $P_0 \to X$. Let $S_0
\subseteq P_0$ be $P_0 \smallsetminus Q_0$ (with the reduced closed substack
structure).  Call $P$ the fiber product
$P_0^t = P_0 \times_{X}
\dots
\times_X P_0$ of $t$ copies of $P_0$. Then the projections $P \to P_0$ are
representable, hence the complement of $S_0^t = S_0 \times_{X} \dots
\times_X S_0$ is an algebraic space. The dimension of $S_0^t$ is
at most equal to $t(d-1) + \dim X$, while the fiber dimension of $P \to X$
is $td$; so for large $t$ we get a smooth projective morphism $P \to X$
with a representable open substack $Q \subseteq P$ such that the
complement of $Q$ in $P$ has dimension less than the fiber dimension
of $\pi$.

Let $U$ be the moduli space of $P$ and let $V$ be the moduli
space of $X$. We have an induced morphism $f\colon U\to V$ which is proper.
By hypothesis, $V$ is quasi-projective.
Let us show that $U$ is quasi-projective as well. This can be shown using
Geometric Invariant Theory, or with the following more elementary method.

\begin{lemma} \label{lem:lbdescend}
Let $\mathcal{L}$ be an invertible sheaf on a
Deligne--Mumford stack $T$,
separated and of finite type over a noetherian base scheme,
and let $U$ be the moduli space of $T$.
Then some power
$\mathcal{L}^{\otimes d}$ (with $d>0$)
is the the pullback of an invertible sheaf $\mathcal{M}$ on $U$.
\end{lemma}

\begin{proof}
First of all, the statement is equivalent to saying that, if we call $\pi
\colon T \to U$ the canonical homomorphism, $\pi_*(\mathcal{L}^d)$ is an
invertible sheaf on $U$, and the adjunction map
$\pi^*\pi_* (\mathcal{L}^d) \to \mathcal{L}^d$ is an isomorphism.
This is a local question in the \'etale topology, and by
\cite[Proposition~4.2]{KM}
we may
assume that $T$ is of the form $[S/G]$, where $G$ is a finite group acting
on an affine scheme $S$. There is a spectral sequence
   \[
   E_2^{pq} = \mathrm{H}^p\bigl(G, \mathrm{H}^q(S,
   \mathbb{G}_{\mathrm{m}})\bigr)
   \Longrightarrow \mathrm{H}^{p+q}([S/G], \mathbb{G}_{\mathrm{m}})
   \]
showing that the kernel of the pullback on Picard groups $\pic [S/G] \to
\pic S$ is $\mathrm{H}^1 \bigl(G, \mathcal{O}^*(S)\bigr)$. By shrinking $S$
we may assume that the pullback of $\mathcal{L}$ to $S$ is trivial; this
means that $\mathcal{L}$ comes from $\mathrm{H}^1 \bigl(G,
\mathcal{O}^*(S)\bigr)$, so some tensor power is trivial.
But $\pi_* \mathcal{O}_T = \mathcal{O}_U$, by definition of moduli space,
so this proves the lemma.
\end{proof}

By Lemma \ref{lem:lbdescend},
there exists an invertible sheaf $\mathcal{M}$ on $U$
whose pullback to $P$ is ample relative to $P \to X$.  We claim that
$\mathcal{M}$ is ample relative to $U \to V$. Again, this is a local
question in the \'etale topology on $V$, so we may assume that $X$ is of
the form $[S/G]$, where $G$ is a finite group and $S$ an affine scheme. We
set $T = S \times_{X} P$, so that $P = [T/G]$. The pullback of
$\mathcal{M}$ to $T$ is ample, and the projection $T \to T/G = U$ is finite
and surjective; hence $\mathcal{M}$ is ample on $U$, and $U$ is
quasiprojective.

By repeated applications of Lemma \ref{lem.bertini}, there is a
complete intersection subscheme $Z\subset U$ with $Z\to V$ surjective
and finite; by dimension reasoning we may take $Z$ to be disjoint from the
image in $U$ of $P\smallsetminus Q$.
Since $Q$ is representable, the morphism $P\to U$ restricts to
an isomorphism of $Q$ to its image, hence $Z$ lifts to a
(representable) substack of $P$, also a complete intersection.
Since $P\to X$ is smooth, by the local criterion for flatness it follows
that $Z$ is flat over $X$.
Also, the singular locus of $Z$ is the pre-image of the singular locus of
$X$.
\end{proof}

\begin{proof}[Proof of Theorem~\ref{thm.tfae}]
First we show that (a) implies (b).
Suppose $\beta\in \mathrm{H}^2(X,\G_m)$ is $n$-torsion, and we want to show
$\beta$ is in the image of the Brauer map.  It suffices to show this
after finite flat pullback, so in the case $\chara k=p$ with $p$ dividing $n$
we take $\varphi\colon X\to X$ to be a suitable iteration of the Frobenius map
and by considering $\varphi^*\beta$ we are reduced to the case
$n$ is relatively prime to $p$.

Now the associated gerbe banded by the $n{\rm th}$ roots of unity is
a tame Deligne--Mumford stack with moduli space $X$, hence is a quotient
stack by (a).
Then by \cite[Theorem~3.6]{EHKV}, $\beta$ lies in the image of the Brauer
map.

Now we show (b) implies (a).
Let $X$ be as in statement (a).
Then, first, $X$ is a gerbe over a smooth separated stack $Y$
such that $Y$ has trivial generic stabilizer and thus by
\cite[Theorem~2.18]{EHKV}, $Y$ is a quotient stack.
More precisely, if $I\to X$ denotes the inertia stack and $U$ is the
open substack of $X$ where the morphism $I\to X$ is flat, with inverse image
$I_U$, then the closure $J$ of $I_U$ in $I$ is
\'etale over $F$, and
$Y$ is the rigidification of
$X$ along $J$ as in \cite[Section~5.1]{ACV}.
The stacks $X$ and $Y$ have the same moduli space.

By Theorem \ref{thm.fifl},
there is a smooth scheme $Z$ and a finite flat surjective morphism
$Z\to Y$, with $Z$ smooth.
The fiber product $W:=X\times_YZ$ is a tame gerbe over $Z$. We claim that,
after replacing
$Z$ by a finite \'etale cover, we can produce a
finite representable \'etale cover $W'\to W$, so that
$W'$ is
a gerbe banded by a product of groups of roots of unity over $Z$.

This can be checked as follows. First of all, we may assume that $Z$ is
connected; then the group of automorphisms of any two geometric points are
isomorphic. Let $G$ be such a group; $X$ is a $G$-gerbe. Then
(\cite[Proposition~3.5]{EHKV}) we may replace $Z$ by a finite \'etale
cover and assume that $W \to Z$ is banded by $G$. After a further cover,
we also assume that the center
$H$ of $G$ is isomorphic to a product of groups of roots of unity over
$Z$. Take $W'$ to be the stack of equivalences of the trivial $G$-gerbe
with $W$ which induce the identity map on bands;
by \cite[Section~1.2.5]{DD} (see also \cite{Gi}) the stack W' is banded by $H$.
There is a tautological evaluation map $W'\to W$,
and this is finite, \'etale, and representable.

This $W'$ is a fiber product $W_1 \times_Z \dots \times_Z
W_r$, where each $W_i$ is a gerbe over $Z$ banded by a group of roots of
unity. Then, by the hypothesis (b) and the main result of \cite{EHKV}, each
$W_i$ is a quotient stack, hence so is $W'$;
by  \cite[Lemma~2.13]{EHKV},
$X$ is a quotient stack as well.
\end{proof}

\begin{proof}[Proof of Theorem~\ref{thm.affinemod}]
The argument is just as in the proof that (b) implies (a) in
Theorem~\ref{thm.tfae}, except that we invoke surjectivity of the
Brauer map for affine schemes, or schemes that are unions of two affines
(\cite{Ga}) in place of the hypothesis (b).
\end{proof}


\begin{thebibliography}{12}

\bibitem[1]{ACV} D. Abramovich, A. Corti, and A. Vistoli, Twisted
bundles and admissible covers, math.AG/0106211, \emph{Comm. Algebra}, to
appear.

\bibitem[2]{AV} D. Abramovich, A. Vistoli, Compactifying the space
of stable maps, \emph{J. Amer. Math. Soc.} \textbf{15} (2002), 27--75.

\bibitem[3]{Co} J.-L. Colliot-Th\'el\`ene, Birational invariants, purity
and the Gersten conjecture, in \emph{$K$-theory and Algebraic Geometry:
Connections with quadratic forms and division algebras (Santa Barbara,
1992)}, \emph{Proc. Sympos. Pure Math.} {\bf 58}, Part 1, Amer. Math. Soc.,
Providence, 1995, 1--64.

\bibitem[4]{DD} P. Debes and J.-C. Douai, Gerbes and covers,
\emph{Comm. Algebra} \textbf{27} (1999), 577--594.


\bibitem[5]{EG} D. Edidin and W. Graham, Equivariant intersection theory,
\emph{Invent. Math.} {\bf 131} (1998), 595--634.

\bibitem[6]{EHKV} D. Edidin, B. Hassett, A. Kresch, and A. Vistoli,
Brauer groups and quotient stacks,
\emph{Amer. J. Math.} {\bf 123} (2001), 761--777.

\bibitem[7]{Ga} O. Gabber, Some theorems on Azumaya algebras, in
\emph{The Brauer Group (Les Plans-sur-Bex, 1980)},
\emph{Lecture Notes in Math.} {\bf 844,} Springer-Verlag, Berlin, 1981,
129--209.

\bibitem[8]{Gi} J. Giraud, \emph{Cohomologie non
ab\'elienne}, Springer-Verlag, Berlin, 1971.

\bibitem[9]{Gr} A. Grothendieck, Le groupe de Brauer I, II, III,
\emph{Dix Expos\'es sur la Cohomologie des Sch\'emas},
\emph{Adv. Stud. Pure Math.} {\bf 3}, North-Holland, Amsterdam, 1968,
46--188.

\bibitem[10]{Lo}
E. Looijenga, Smooth Deligne--Mumford compactifications by means of Prym
level structures, \emph{J. Algebraic Geom.} {\bf 3} (1994), 283--293.

\bibitem[11]{KM} S. Keel and S. Mori, Quotients by groupoids,
\emph{Ann. of Math. (2)} {\bf 145} (1997), 193--213.

\bibitem[12]{Mi} J. S. Milne, \emph{Etale Cohomology},
Princeton Univ. Press, Princeton, 1980.

\bibitem[13]{Mu} D. Mumford, Towards an enumerative geometry of the moduli
space of curves, in \emph{Arithmetic and Geometry},
\emph{Progr. Math.}, vol. 36, Birkh\"auser, Boston, 1983, 271--328.

\bibitem[14]{PDJ} M. Pikaart and A. J. de Jong, Moduli of curves with
non-abelian level structure, in \emph{The moduli space of curves (Texel
Island, 1994)}, \emph{Progr. Math.} {\bf 129}, Birkh\"auser, Boston, 1995,
483--509.

\bibitem[15]{To} B. Totaro, The resolution property for schemes and stacks,
\emph{J. Reine Angew. Math.}, to appear.
\end{thebibliography}
\end{document}